# A STUDY ON CLIQUE INVARIANTS OF JACO-TYPE GRAPHS

M. Jerlin Seles

Research Scholar, Department of Mathematics

Nirmala College for Women

Coimbatore, India

jerlinseles@gmail.com

U. Mary

Department of Mathematics

Nirmala College for Women

Coimbatore, India

marycbe@gmail.com

Johan Kok

Tshwane Metropolitan Police Department

City of Tshwane, Republic of South Africa

Kokkiek2@tshwane.gov.za

**Abstract:**

The first study related to this concept was that on primitive holes [2]. A primitive hole is indeed the study of the existence of 3-cliques. In this paper we report on research in respect of clique parameters and related properties thereof for certain Jaco-type graphs.

## 1. INTRODUCTION

For basic definitions and results not mentioned in this paper, we refer to [1,2,3,5]. A Jaco graph $J_n(f(1))$ as defined in the earlier work is exactly a linear Jaco graph $J_n(f(x))$ with $f(x) = x$. (See [1,3]). Also the definition of a linear Jaco graph (all Jaco graphs for that matter) incorporates an integer valued linear function. The function value of the vertex subscript say, $f(i)$ of vertex $v_i$ determines the total vertex degree $d(v_i)$.

Completely separate (in definition) from linear Jaco graphs we have Jaco-type graphs. These are the graphs for which a non-negative integer sequence defines only the out-degree of a vertex. So entry $a_i$ of a non-negative integer sequence is the out-degree, $d^+(v_i)$ of vertex $v_i$.

The graphs generally show diagrammatical resemblance to linear Jaco graphs when sketched and many similar results hold between them. In fact there is a non-empty intersection between the sets of these graphs. It means some (perhaps all, but not proven yet) linear Jaco graphs are Jaco-type graphs. It was found that the linear Jaco graph $J_n(x)$ (remember now $f(x) = x$) is also a Jaco-type graph and the sequence defining the out-degrees is found in paper [5].

For many Jaco-type graphs one cannot find a linear function defining the total vertex degree. For example, a linear function defining the total vertex degrees of the Fibonacci Jaco-type graph is not yet found. So until proven otherwise the Fibonacci Jaco-type graph is considered not to be a linear Jaco graph.

## 2. JACO-TYPE GRAPHS AND CLIQUE PARAMETERS

### 2.1 Basic Results for Certain Jaco-type Graphs

The linear Jaco graph for $f(x) = x$ is indeed a Jaco-type graph on the sequence, $\{a_n\} = n - \left\lfloor \frac{2(n+1)}{3+\sqrt{5}} \right\rfloor, n = 1,2,3,...$ A closely related Jaco-type graph is that on the positive integer sequence $s_1 = \{a_n\} = 1,2,3,...$ and for brevity be denoted, $J_\infty(s_1)$ or $J_n(s_1)$ for the infinite and finite Jaco-type graphs respectively. Figures 2.1.1 and 2.1.2 below depict the Jaco-type graphs for $J_8(s_1)$ and $J_{12}(s_2)$ respectively, with $s_2 = \{f_i\}$, for $f_0 = 0, f_1 = 1, f_2 = 1, f_3 = 2,...$ The aforesaid sequence is the well-known Fibonacci sequence. (Also see [5]).

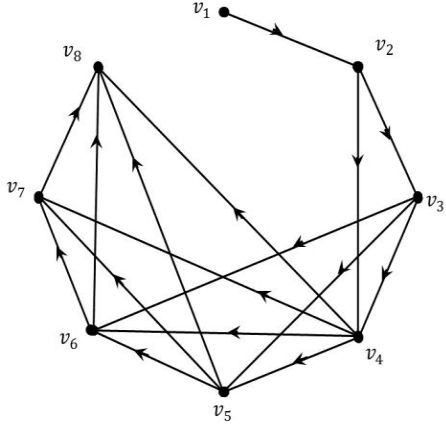 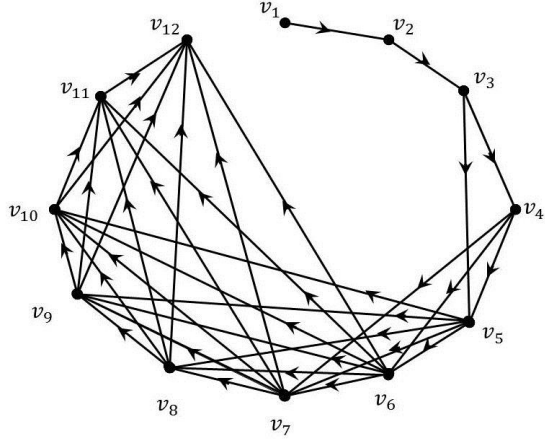

Figure 2.1.1 [5]  Figure 2.1.2 [5]

The underlying Jaco-type will be denoted $J_n^*(s_k)$ and when the context is clear we refer to both graphs as Jaco-type graphs. Also the use of terminology such as arc versus edge will be understood to refer to the directed Jaco-type graph versus its underlying graph, respectively.

The Jaco-type graph $J_8(s_1)$ has $\Delta(J_8(s_1)) = 6$ and $J(J_8(s_1)) = \{v_4\}$. It has girth 3 and circumference 5.

Similarly, the Jaco-type graph $J_{12}(s_2)$ has $\Delta(J_{12}(s_2)) = 8$ and $J(J_{12}(s_2)) = \{v_6, v_7\})$. It has girth 3 and circumference 7.

### 2.1.1 Basic properties and results of Jaco-type graphs

A comprehensive study of the properties of $J_n(s_1)$, is found in [5]. Properties for the Jaco-type graph are given without proof in [5] because they can easily be verified against the definition of infinite Jaco-type graph, $J_\infty(\{a_n\})$.

**Lemma 2.1.1**

The in-degree $d^-(v_i)$, for any vertex $v_i$ found in both $J_\infty(s_k)$ and $J_n(s_k)$ remains a constant for any given integer sequence, $\{s_k\}$.

**Proof:**

Consider $J_n(s_k)$ for any integer sequence $\{s_k\}$ and any integer $n \in N$. Also consider any vertex $v_i$ with $d^-(v_i) = l > 0$. Now extend to the Jaco-type graph $J_{n+1}(s_k)$. Clearly, $d^-(v_{n+1}) = t > 0$ and $d^+(v_{n+1}) = 0$. Hence, only out-arcs were added to some corresponding vertices, possibly $v_i$ as well, but $v_i$ did not get an additional in-arc. Therefore, in the Jaco-type graph $J_{n+1}(s_k)$ we have $d^-(v_i) = l$. Hence, through immediate induction it follows that $d^-(v_i) = l$, a constant in all Jaco-type graphs corresponding to the integer sequence, $\{s_k\}$. □

**Proposition 2.1.2**

A Jaco-type graph $J_n(s_k)$ which, for the smallest $i$ has a vertex $v_i$ with $d^+(v_i) > 1$ and $n \geq i+2$, has girth, $g(J_n(s_k)) = 3$.

**Proof:**

It follows from definition of $J_\infty(\{a_n\})$ that a Jaco-type graph $J_n(s_k)$ which, for the smallest $i$ has a vertex $v_i$ with $d^+(v_i) > 1$ and $n = i+2$, has exactly one smallest cycle $C_3$. This cycle will remain the shortest cycle in all Jaco-type graphs, $J_n(s_k), n \geq i+2$. Hence, the result. □

**Proposition 2.1.3**

For any given integer sequence $\{s_k\}$ the circumference of the Jaco-type graph $J_n(s_k)$ which, for the smallest $i$ has a vertex $v_i$ with $d^+(v_i) > 1$ and $n \geq i+2$, is equal to the clique number $\omega(J_n(s_k))$.

**Proof:**

In general, the clique number equals the order of a maximum clique which in itself is maximal, such maximum clique *per se* has circumference of maximum cycle length. Because the largest hole (chordless cycle) of a Jaco-type graph, if it exists, is $C_3$ it follows that the circumference of such maximum clique is equal to the circumference of its supergraph, $J_n(s_k)$. □

**Proposition 2.1.4**

For any given non-decreasing integer sequence $s_k$ the Jaco-type graph $J_n(s_k)$ which, for the smallest $i$ has the arc $(v_i, v_n)$, has clique cover number, $c(J_n(s_k)) = i$.

**Proof:**

From the definition of infinite Jaco-type graph, $J_\infty(\{a_n\})$ it follows that the vertex set $\{v_i, v_{i+1}, v_{i+2}, ..., v_n\}$ induces a maximal clique. Similarly the vertex sets,

$\{v_{i-1}, v_i, v_{i+1}, ..., v_{(i-1)+a_{(i-1)}}\}, \{v_{i-2}, v_{i-1}, v_i, ..., v_{(i-2)+a_{(i-2)}}\}, ..., \{v_1, v_2, v_3, ..., v_{1+a_1}\}$ all induce maximal cliques. Clearly the union of the corresponding clique vertex sets covers $V(J_n(s_k))$. Therefore, $c(J_n(s_k)) \leq i$.

Assume, $c(J_n(s_k)) < i$, It implies that at least two of the cliques must form a single clique which implies not all cliques were maximal. The aforesaid is a contradiction hence, $c(J_n(s_k)) = i$. □

## 2.2 On *l*-Cliques in Jaco-type Graphs

It is noted from figures 2.1.1 and 2.1.2 that Jaco-type graphs are typically a structural combination of numerous complete graphs or put differently, of numerous *k*-cliques. Therefore, a good understanding of clique parameters related to complete graphs is important before analyzing same for Jaco-type graphs.

From the definition of a complete graph it follows immediately that the induced subgraph $\langle X \rangle, X \subseteq V(K_n)$ is a complete graph as well. Denote the number of distinct $l$-cliques imbedded in $K_n$ by $\eta^{K_l}(K_n)$.

Proposition 2.2.1 shows that the result stated in *WolframMathWorld.com* [8] is incorrect.

**Proposition 2.2.1**

The total number of distinct $l$-cliques $l = 1,2,3,...,n$ imbedded in a complete graph $K_n$ is

$$\sum_{l=1}^{n} \binom{n}{l}.$$

**Proof:**

The number of ways to select $l$-subsets (cardinality $l$) for a set with cardinality $n$ is given by $\binom{n}{l}$. Because $l$-cliques, $l = 1,2,3,...,n$ exist in $K_n$ the result follows immediately. □

**Definition 2.2.2**

The **join** of graphs $G$ and $H$ is the graph obtained by adding edges to link each vertex of $G$ to all vertices of $H$.

**Proposition 2.2.3 (Generalized Clique Theorem)**

For any finite graph of order, $n \geq 1$, we have $\eta^{K_{l+1}}(G + K_1) = \eta^{K_{l+1}}(G) + \eta^{K_l}(G), 0 \leq l \leq n.$

**Proof:**

(i) Because all graphs have an empty clique, $\eta^{K_0}(G) = 1$.

Therefore, $\eta^{K_1}(G + K_1) = \eta^{K_1}(G) + 1 = n + 1$.

(ii) Consider any clique $C_i = K_l$, $0 \leq l \leq n$ and $1 \leq i \leq \eta^{K_l}(G)$ of, $G$.

Clearly $C_i + K_1 = K_{l+1}$. Hence in $G + K_1$, exactly $\eta^{K_l}(G)$ new cliques of order $l+1$ will be created. Since, the graph $G$ itself has $\eta^{K_{l+1}}(G) \geq 0$ such cliques we have that $\eta^{K_{l+1}}(G+K_1) = \eta^{K_{l+1}}(G) + \eta^{K_l}(G)$. The general result follows through immediate induction. □

Note that a complete graph of order $n$ is the special case: $K_n = ((((K_1+K_1)+K_1)+K_1)+...+K_1)$, a $(n-1)$-fold join. Applying Proposition 2.2.3 to a complete graph, $K_{n+1}$, $n \geq 0$ offers a complete graph specific proof. We present it as a corollary together with alternative proof technique.

**Corollary 2.2.4**

For a complete graph $K_n$, $n \geq 1$ we have, $\eta^{K_{l+1}}(K_n) = \eta^{K_{l+1}}(K_{n-1}) + \eta^{K_l}(K_{n-1})$ with $0 \leq l \leq n$.

**Proof:**

Consider any complete graph $K_m$, $m \geq 1$. Now consider $\binom{m}{l} + \binom{m}{l-1}$.

$$\binom{m}{l} + \binom{m}{l-1} = \frac{m!}{l!(m-l)!} + \frac{m!}{(l-1)!(m-l+1)!}$$

$$= \frac{m!(m-l+1)}{l!(m-l+1)(m-l)!} + \frac{m!}{l(l-1)!(m-l+1)!}$$

$$= \binom{m+1}{l}.$$

Hence, through immediate induction we have, $\eta^{K_{l+1}}(K_n) = \eta^{K_{l+1}}(K_{n-1}) + \eta^{K_l}(K_{n-1}), n \geq 1$ with $0 \leq l \leq n$. □

**Theorem 2.2.5**

For any complete graph $K_n$, $n \geq 1$ we have $\eta^{K_j}(K_n) = \eta^{K_{n-j}}(K_n)$, $1 \leq j \leq n$ (inclusive of the empty-clique).

**Proof:**

By the Proposition 2.2.2, $\eta^{K_j}(K_n) = \binom{n}{j}$ for, $1 \leq j \leq n$.

We have, $\binom{n}{j} = \dfrac{n!}{j!(n-j)!} = \dfrac{n!}{(n-j)!(n-(n-j))!} = \binom{n}{n-j}$,

and since $\eta^{K_{n-j}}(K_n) = \binom{n}{n-j}$ by definition, the result that $\eta^{K_j}(K_n) = \eta^{K_{n-j}}(K_n)$ for $1 \leq j \leq n$ follows. $\square$

Table 1 below depicts the number of cliques (excluding the empty clique) of complete graphs $K_n, 1 \leq n \leq 10$ and $1 \leq l \leq 10$. Important matrix properties will follow table 1.

Table 1

| $K_n$, $n=$ | $\eta^{K_1}$ | $\eta^{K_2}$ | $\eta^{K_3}$ | $\eta^{K_4}$ | $\eta^{K_5}$ | $\eta^{K_6}$ | $\eta^{K_7}$ | $\eta^{K_8}$ | $\eta^{K_9}$ | $\eta^{K_{10}}$ |
|---|---|---|---|---|---|---|---|---|---|---|
| 1 | 1 | 0 | 0 | 0 | 0 | 0 | 0 | 0 | 0 | 0 |
| 2 | 2 | 1 | 0 | 0 | 0 | 0 | 0 | 0 | 0 | 0 |
| 3 | 3 | 3 | 1 | 0 | 0 | 0 | 0 | 0 | 0 | 0 |
| 4 | 4 | 6 | 4 | 1 | 0 | 0 | 0 | 0 | 0 | 0 |
| 5 | 5 | 10 | 10 | 5 | 1 | 0 | 0 | 0 | 0 | 0 |
| 6 | 6 | 15 | 20 | 15 | 6 | 1 | 0 | 0 | 0 | 0 |
| 7 | 7 | 21 | 35 | 35 | 21 | 7 | 1 | 0 | 0 | 0 |
| 8 | 8 | 28 | 56 | 70 | 56 | 28 | 8 | 1 | 0 | 0 |
| 9 | 9 | 36 | 84 | 126 | 126 | 84 | 36 | 9 | 1 | 0 |
| 10 | 10 | 45 | 120 | 210 | 252 | 210 | 120 | 45 | 10 | 1 |

The entries of the table 1 can be written in matrix form as:

$$A = \begin{pmatrix} 1 & 0 & 0 & 0 & \cdots & 0 \\ 2 & 1 & 0 & 0 & \cdots & 0 \\ 3 & 3 & 1 & 0 & \cdots & 0 \\ 4 & 6 & 4 & 1 & \cdots & 0 \\ \cdots & \cdots & \cdots & \cdots & \cdots & \cdots \\ n & a_{n-1,1}+a_{n-1,2} & a_{n-1,2}+a_{n-1,3} & \cdots & \cdots & 1 \end{pmatrix}.$$

Note that the first column has the entries, $a_{i,1} = i$, $i = 1,2,3,\ldots,n$. All other entries $a_{i,j} = a_{i-1,j-1} + a_{i-1,j}$. Note that matrix $A$ is a lower triangular matrix. Since the determinant of a lower triangular matrix is the product of the diagonal elements of the matrix we have: $\det(A) = 1$. Because of Corollary 2.2.4, inverse matrix, $A^{-1}$ follows remarkably easily.

**Proposition 2.2.6**

The inverse of matrix A is given by:

$$A^{-1} = \begin{pmatrix} +1 & 0 & 0 & 0 & \cdots & 0 \\ -2 & +1 & 0 & 0 & \cdots & 0 \\ +3 & -3 & +1 & 0 & \cdots & 0 \\ -4 & +6 & -4 & +1 & \cdots & 0 \\ \cdots & \cdots & \cdots & \cdots & \cdots & \cdots \\ +n & -(a_{n-1,1}+a_{n-1,2}) & +(a_{n-1,2}+a_{n-1,3}) & \cdots & \cdots & +1 \end{pmatrix} \text{ if } n \text{ is odd,}$$

and

$$A^{-1} = \begin{pmatrix} +1 & 0 & 0 & 0 & \cdots & 0 \\ -2 & +1 & 0 & 0 & \cdots & 0 \\ +3 & -3 & +1 & 0 & \cdots & 0 \\ -4 & +6 & -4 & +1 & \cdots & 0 \\ \cdots & \cdots & \cdots & \cdots & \cdots & \cdots \\ -n & +(a_{n-1,1}+a_{n-1,2}) & -(a_{n-1,2}+a_{n-1,3}) & \cdots & \cdots & +1 \end{pmatrix} \text{ if } n \text{ is even.}$$

**Proof:**

All entries of matrix A has the form $a_{i,j} = a_{i-1,j-1} + a_{i-1,j}$. The validity of the entry form follows directly from Proposition 2.2.2 applied to the complete graph, $K_n$, $n \geq 1$.

**Case 1:** If $n$ is odd, let:

$$A^{-1} = \begin{pmatrix} +1 & 0 & 0 & 0 & \cdots & 0 \\ -2 & +1 & 0 & 0 & \cdots & 0 \\ +3 & -3 & +1 & 0 & \cdots & 0 \\ -4 & +6 & -4 & +1 & \cdots & 0 \\ \cdots & \cdots & \cdots & \cdots & \cdots & 0 \\ +n & -(a_{n-1,1} + a_{n-1,2}) & +(a_{n-1,2} + a_{n-1,3}) & \cdots & \cdots & +1 \end{pmatrix}.$$

(i) The first $n$-row $(a_{1,j}), 1 \leq j \leq n$ in matrix A has the form $r_1 = (1,0,0,...,0)$ and all other $n$-row in matrix A has the form $r_i = (i, a_{i,j}), 2 \leq i \leq n, 2 \leq j \leq n$ with $a_{i,j} = a_{(i-1),(j-1)} + a_{(i-1),j}$. Every column $c_{i,j}, 1 \leq i \leq n, 1 \leq j \leq n$ of matrix $A^{-1}$ has the entries of the columns of matrix $A$ together with alternating + and - sign beginning from each diagonal entry +1. Clearly, $r_{i,j}.c_{i,j} = 1$ and $r_{i,j}.c_{k,j} = 0, i \neq k$. Therefore, $A.A^{-1} = I$.

(ii) The first $n$-row $(a_{1,j}), 1 \leq j \leq n$ in matrix $A^{-1}$ has the form $r_1 = (1,0,0,...,0)$ and all other $n$-rows in matrix $A^{-1}$ has the form $r_i = (\pm i, a_{i,j}), 2 \leq i \leq n, 2 \leq j \leq n$ and $+i$ if and only if $i$ is odd with $|a_{i,j}| = |a_{(i-1),(j-1)}| + |a_{(i-1),j}|$. Clearly, $r_{i,j}.c_{i,j} = 1$ and $r_{i,j}.c_{k,j} = 0, i \neq k$. Therefore, $A^{-1}.A = I$.

**Case 2:** If $n$ is even, let

$$A^{-1} = \begin{pmatrix} +1 & 0 & 0 & 0 & \cdots & 0 \\ -2 & +1 & 0 & 0 & \cdots & 0 \\ +3 & -3 & +1 & 0 & \cdots & 0 \\ -4 & +6 & -4 & +1 & \cdots & 0 \\ \cdots & \cdots & \cdots & \cdots & \cdots & 0 \\ -n & +(a_{n-1,1}+a_{n-1,2}) & -(a_{n-1,2}+a_{n-1,3}) & \cdots & \cdots & +1 \end{pmatrix} \text{ if } n \text{ is even.}$$

This case follows similar to Case 1. Therefore, $A.A^{-1} = A^{-1}.A = I$. □

## 2.3 Vertex Clique Degrees of Certain Jaco-type Graphs

We begin this section by stating certain important results in respect of vertex clique degrees for complete graphs. This is followed by applications to the sequence of positive integers, the Fibonacci Jaco-type graph, modulo $k$ Jaco-type graph and the set Jaco-type graph. We begin with an important theorem.

**Theorem 2.3.1**

For a complete graph $K_n$, $n \geq 1$ we have that $d^{K_l}(v_i) = \dfrac{l.\eta^{K_l}(K_n)}{n}$.

**Proof:**

For the complete graph $K_1$ we have $d^{K_1}(v_1) = \dfrac{1.1}{1} = 1$. For the complete graph $K_2$ we have that $d^{K_1}(v_i) = \dfrac{1.2}{2} = 1, i = 1,2$ and $d^{K_2}(v_i) = \dfrac{2.1}{2} = 1, i = 1,2$. Hence, the result holds for $n = 1,2$.

Assume the result holds for $1 \leq l \leq m$. So, $d^{K_l}(v_i) = \dfrac{l.\eta^{K_l}(K_n)}{n}, 1 \leq l \leq m$.

Consider the complete graph $K_{m+1}$. Clearly, $d^{K_1}(v_i) = \dfrac{1.(m+1)}{m+1} = 1, 1 \leq i \leq m+1$. By the definition of a complete graph all vertex degrees are equal.

(i.e.)., $d^{K_2}(v_i) = m, 1 \leq i \leq m+1$.

Also, $d^{K_2}(v_i) = \dfrac{2.\eta^{K_2}(K_{m+1})}{m+1} = \dfrac{2 \cdot \frac{1}{2}(m+1)m}{m+1} = m.$

Therefore, $d^{K_2}(v_i) = \dfrac{2.\eta^{K_2}(K_{m+1})}{m+1}, 1 \leq i \leq m+1.$ Thus the result holds for the complete graph $K_{m+1}$ in respect of $d^{K_1}(v_i)$ and $d^{K_2}(v_i), 1 \leq i \leq m+1$.

Now consider any $l \leq m$. We have that vertex $v_{m+1}$ induces $\binom{m}{l-1}$ complete graphs $K_l$.

Hence, $d^{K_l}(v_i) = \dfrac{m!}{(l-1)!(m-l+1)!}, 1 \leq i \leq m+1.$

Also,

$$\dfrac{l.\eta^{K_l}(K_{m+1})}{m+1} = \dfrac{l\binom{m+1}{l}}{m+1} = \dfrac{l(m+1)!}{(m+1)l!(m+1-l)!} = \dfrac{m!}{(l-1)!(m-l+1)!} = d^{K_l}(v_i), 1 \leq i \leq m+1.$$

Finally, $d^{K_{m+1}}(v_i) = 1, 1 \leq i \leq m+1$ in $K_{m+1}$.

Therefore, through imbedded induction the result holds for all complete graphs $K_n, n \in N$. □

Note that applying Theorem 2.3.1 to a complete graph $K_n, n \geq 1$ in respect of edges (2-cliques) re-establishes the well-known result that the degree of each vertex is $d^{K_2}(v_i) = \dfrac{2.\eta^{K_2}(K_n)}{n} = \dfrac{1}{2}n(n-1).$ In fact the inverse formula may serve as an generalization to obtain the number of edges of a graph on $n$ vertices is given by $\dfrac{1}{2}\sum_{i=1}^{n}\deg(v_i).$ □

Table 2 depicts the vertex clique degrees corresponding to Table 3.2.1 Note that for a given $K_n$, symmetry ensures that all vertices $v \in V(K_n)$ have equal vertex clique degree in respect of a specific clique size.

Table 2 Vertex Clique Degrees

| $K_n$, $n=$ | $d^{K_1}(v)$ | $d^{K_2}(v)$ | $d^{K_3}(v)$ | $d^{K_4}(v)$ | $d^{K_5}(v)$ | $d^{K_6}(v)$ | $d^{K_7}(v)$ | $d^{K_8}(v)$ | $d^{K_9}(v)$ | $d^{K_{10}}(v)$ |
|---|---|---|---|---|---|---|---|---|---|---|
| 1 | 1 | 0 | 0 | 0 | 0 | 0 | 0 | 0 | 0 | 0 |
| 2 | 1 | 1 | 0 | 0 | 0 | 0 | 0 | 0 | 0 | 0 |
| 3 | 1 | 2 | 1 | 0 | 0 | 0 | 0 | 0 | 0 | 0 |
| 4 | 1 | 3 | 3 | 1 | 0 | 0 | 0 | 0 | 0 | 0 |
| 5 | 1 | 4 | 6 | 4 | 1 | 0 | 0 | 0 | 0 | 0 |
| 6 | 1 | 5 | 10 | 10 | 5 | 1 | 0 | 0 | 0 | 0 |
| 7 | 1 | 6 | 15 | 20 | 15 | 6 | 1 | 0 | 0 | 0 |
| 8 | 1 | 7 | 21 | 35 | 35 | 21 | 7 | 1 | 0 | 0 |
| 9 | 1 | 8 | 28 | 56 | 70 | 56 | 28 | 8 | 1 | 0 |
| 10 | 1 | 9 | 36 | 84 | 126 | 126 | 84 | 36 | 9 | 1 |

From Table 2 an interesting theorem follows.

**Theorem 2.3.2**

For a complete graph $K_n$ we have: $d^{K_l}(v_i) = \dfrac{\prod_{j=1}^{l-1} n-j}{n!}$, $n = 1,2,3,...$, with $2 \leq l \leq n$ and $1 \leq i \leq n$.

**Proof:**

Clearly by default the clique degree, $d^{K_1}(v_i) = 1$ for all $K_n$, $n = 1,2,3,...$ It follows easily that:

$$d^{K_2}(v_i) = n-1, n \geq 2,$$
$$d^{K_3}(v_i) = \frac{(n-1)(n-2)}{2!}, n \geq 3,$$
$$d^{K_4}(v_i) = \frac{(n-1)(n-2)(n-3)}{3!}, n \geq 4.$$

Assume the result holds for $d^{K_l}(v_i), n \geq l$. So by the induction assumption

$$d^{K_l}(v_i) = \frac{\prod_{j=1}^{l-1} n-j}{n!}, n \geq l.$$ Now consider, $K_{l+1}$, and it follows that:

$$d^{K_2}(v_i) = (l-1)+1 = (l+1)-1,$$
$$d^{K_3}(v_i) = \frac{((l-1)+1)((l-2)+1)}{2!} = \frac{((l+1)-1)((l+1)-2)}{2!},$$
$$d^{K_4}(v_i) = \frac{((l-1)+1)((l-2)+1)((l-3)+1)}{3!} = \frac{((l+1)-1)((l+1)-2)((l+1)-3)}{3!},$$
...
...

$$d^{K_{l+1}}(v_i) = \frac{\prod_{j=1}^{l} (l+1)-j}{(l+1)!}.$$

Through immediate induction it follows that it also holds for the clique degrees of $K_{l+t}, t = 2,3,4,...$

Hence through induction we have for $K_n$, $n \geq 1$, $\forall v_i \in V(K_n)$ $n = 1,2,3,...$ that:

$$d^{K_l}(v_i) = \frac{\prod_{j=1}^{l-1} n-j}{n!}. \qquad \square$$

**Proposition 2.3.3**

For a complete graph, $K_n, n \geq 1$ the maximum clique degree is $d^{K_t}(v_i)$ with:

$$t = \begin{cases} \left\lceil \dfrac{n}{2} \right\rceil, & \text{if } n \text{ is odd}, \\ \dfrac{n}{2} \text{ or } \dfrac{n}{2}+1, & \text{if } n \text{ is even}. \end{cases}$$

**Proof:**

From the definition of the entries of matrix $A$ in the proof of Proposition 2.2.6 it follows immediately that the maximum number of cliques for a complete graph $K_n$ is given by $\eta^{K_t}(K_n)$ with:

$$t = \begin{cases} \left\lceil \dfrac{n}{2} \right\rceil, & \text{if } n \text{ is odd}, \\ \dfrac{n}{2} \text{ or } \dfrac{n}{2}+1, & \text{if } n \text{ is even.} \end{cases}$$

Hence, the result follows directly from Theorem 2.3.1. □

### 2.3.1 Application to the finite positive integer sequence Jaco-type Graph: $J_n(s_1)$ [05]

The infinite Jaco-type graph $J_\infty(s_1)$ is the graph with vertex set $V(J_\infty(s_1)) = \{v_i : i \in N\}$ and the arc set $A(J_\infty(s_1)) = \{(v_i, v_j) : i, j \in N, i < j\}$ such that $(v_i, v_j) \in A(J_\infty(s_1))$ if and only if $2i \geq j$. Note that a finite Jaco-type graph $J_n(s_1)$ in this family is obtained from $J_\infty(s_1)$ by lobbing off all vertices $v_k$, $\forall k > n$ (with incident arcs).

**Theorem 2.3.4**

The finite integer Jaco-type graph $J_n(s_1)$ can be decomposed in exactly the number of maximal cliques $K_2, K_3, ..., K_{\left(\frac{n}{2}\right)}, n \geq 2$ and even and into $K_2, K_3, ..., K_{\left\lfloor \frac{n}{2} \right\rfloor +1}, K_{\left\lfloor \frac{n}{2} \right\rfloor +1}, n \geq 1$ and odd corresponding to the clique cover number $c(J_n(s_1)) = \dfrac{n}{2}$ or $\dfrac{n+1}{2}$.

**Proof:**

For $J_1(s_1)$ we have $K_1$. For $J_2(s_1)$ we have the path which is the one maximal clique $P_2$. For $J_3(s_1)$ we have the path $P_3$ therefore the decomposition into maximal cliques corresponding to the clique cover number is the two cliques $K_2, K_2$ on vertices $v_1, v_2$ and $v_2, v_3$. For $J_4(s_1)$ the decomposition into maximal cliques corresponds to $K_2, K_3$. For $J_5(s_1)$ the decomposition into maximal cliques corresponds to $K_2, K_3, K_3$. For $J_6(s_1)$ the decomposition into maximal cliques corresponds to $K_2, K_3, K_4$. Clearly the result holds for $J_l(s_1), 1 \leq l \leq 6$. Assume the result holds for $1 \leq l \leq m$ and without loss of generality assume that the prime Jaconian vertex say, $v_p$ has defined maximum out-degree $d^+(v_p) = p$ and the arc $(v_p, v_m)$ exists.

Now extend to $J_{m+1}(s_1)$. Clearly, a new maximal clique is created.

If $m$ is even then $c(J_m(s_1)) + 1 = \frac{m}{2} + 1 = \frac{(m+1)+1}{2} = c(J_{m+1}(s_1))$.

By similar reasoning the case for, $m$ is odd, follows. Hence, through induction the result follows. □

**Corollary 2.3.5**

For $J_n(s_1)$, $n$ odd and for $J_{n+1}(s_1)$, the number of maximal cliques in the decomposition is equal.

**Proof:**

The result is an immediate consequence of the immediate induction in the proof of Theorem 2.3.4. □

By applying Theorem 2.3.4 it is possible to determine the number of cliques of any order in the integer Jaco-type to the maximum order that exists for a given $n$ namely, $\left(\frac{n}{2}\right)$ for even $n$ and $\left\lfloor\frac{n}{2}\right\rfloor + 1$ for odd $n$. Because certain pairs of maximal cliques have non-empty vertex intersection the determining of distinct cliques and vertex clique degrees must discount multiple counting.

**Example 2.3.1**

Figure 2.1.1 depicts $J_8(s_1)$ [5]. The decomposition into maximal cliques results in $\frac{8}{2} = 4$ cliques (i.e.), $K_2$ on vertices $v_1, v_2$; $K_3$ on vertices $v_2, v_3, v_4$; $K_4$ on vertices $v_3, v_4, v_5, v_6$ and $K_5$ on vertices $v_4, v_5, v_6, v_7, v_8$.

We know that $\eta^{K_1}(J_8(s_1)) = 8$. If the counting was through the number of maximal cliques then a double count at vertex $v_2$ (common to $K_2$ and $K_3$) must be discounted by 1. The double count of vertex $v_3$ (common to $K_3$ and $K_4$) and the triple count of vertex $v_4$ (common to $K_3, K_4$ and $K_5$) must be discounted for by 1 and 2 respectively. Also the double count of vertices $v_5, v_6$ (common to $K_4$ and $K_5$) must be discounted for by 1 in each count to yield the correct count.

Similarly, $\eta^{K_2}(J_8(s_1)) = 16$. If the counting was through the number of maximal cliques the initial count would be 20. Then discounting for the double count of arc $(v_3, v_4)$ (common to $K_3$ and $K_4$) and for arcs $(v_4, v_5)$, $(v_4, v_6)$, $(v_5, v_6)$ (common to $K_4$ and $K_5$) renders the correct result.

For the 3-cliques we have $\eta^{K_3}(J_8(s_1)) = 14$. Counting through the distinct maximal cliques amounts to 15 cliques, $K_3$. After discounting the double count of the triangle on the vertices $v_4, v_5, v_6$ the correct $\eta^{K_3}(J_8(s_1)) = 14$ is obtained.

Finally, in respect of $\eta^{K_4}(J_8(s_1))$ and $\eta^{K_5}(J_8(s_1))$ no double count occurs.

**Lemma 2.3.6**

In the root integer Jaco-type graph $J_\infty(s_1)$ a maximal say, $K_l, l \geq 2$, intersects with $K_{l+1}$ in respect of $l-1$ vertices, intersects with $K_{l+2}$ in respect of $l-2$ vertices, and consecutively so on, until the intersection with $K_{l+(l-1)}$, with which it intersects in respect of a single vertex.

**Proof:**

From definition of infinite Jaco-type graph, $J_\infty(\{a_n\})$ it follows that for any vertex $v_l \in V(J_\infty(s_1))$ the arcs $(v_l, v_{l+1}), (v_l, v_{l+2}), (v_l, v_{l+3}),...,(v_l, v_{2l})$ exist. Hence, the result follows immediately. $\square$

In applying Lemma 2.3.6 care must be taken to discount maximal cliques which were lobbed off as well as arcs which were lobbed off to reduce the order of certain maximal cliques.

### 2.3.2 Application to the finite Fibonacci Jaco-type Graph: $J_n(s_2)$ [6]

The infinite Jaco-type graph corresponding to Fibonacci sequence, which is also called the **Fibonaccian Jaco-type graph** $J_\infty(s_2)$ and is defined by the vertex set $V(J_\infty(s_2)) = \{v_i : i \in N\}$ and the arc set $A(J_\infty(s_2)) = \{(v_i, v_j) : i, j \in N, i < j\}$ and $(v_i, v_j) \in A(J_\infty(s_2))$ if and only if $i + f_i \geq j$. Note that a finite Jaco-type graph $J_n(s_2)$ in this family is obtained from $J_\infty(s_2)$ by lobbing off all vertices $v_k$, $\forall k > n$. Figure 2.1.2 depicts $J_{12}(s_2)$. The table below depicts the number of cliques of all cliques sizes found in Fibonaccian Jaco-type graphs, $1 \leq n \leq 12$.

Table3

| $J_n(s_2)$, $n =$ | $\eta^{K_1}$ | $\eta^{K_2}$ | $\eta^{K_3}$ | $\eta^{K_4}$ | $\eta^{K_5}$ | $\eta^{K_6}$ | $\eta^{K_7}$ |
|---|---|---|---|---|---|---|---|
| 1 | 1 | 0 | 0 | 0 | 0 | 0 | 0 |
| 2 | 2 | 1 | 0 | 0 | 0 | 0 | 0 |
| 3 | 3 | 2 | 0 | 0 | 0 | 0 | 0 |
| 4 | 4 | 3 | 0 | 0 | 0 | 0 | 0 |
| 5 | 5 | 5 | 1 | 0 | 0 | 0 | 0 |
| 6 | 6 | 7 | 2 | 0 | 0 | 0 | 0 |
| 7 | 7 | 10 | 5 | 1 | 0 | 0 | 0 |
| 8 | 8 | 13 | 8 | 2 | 0 | 0 | 0 |

| 9 | 9 | 17 | 12 | 6 | 1 | 0 | 0 |
|---|---|---|---|---|---|---|---|
| 10 | 10 | 22 | 22 | 16 | 6 | 1 | 0 |
| 11 | 11 | 27 | 32 | 26 | 11 | 1 | 0 |
| 12 | 12 | 33 | 47 | 46 | 17 | 7 | 1 |

**Lemma 2.3.7**

For the Fibonaccian Jaco-type graph $J_{n+1}(s_2)$ with $d^-(v_{n+1}) = l$ then

$$\eta^{K_i}(J_{n+1}(s_2)) = \binom{n+1}{i} + \eta^{K_i}(J_n(s_2)), \ 2 \leq i \leq l.$$

**Proof:**

Clearly, the number of combinations, $\binom{n+1}{i}$ of in-arcs of vertex $v_{n+1}$ corresponds to the additional cliques, $K_i$ on expanding from $J_n(s_2)$ to $J_{n+1}(s_2)$. Therefore, the result follows as a derivative of Theorem 2.2.3. □

### 2.3.3 Application to the finite Modulo $k$ Jaco-type Graph: $J_n(s_3)$ [6]

The notion of modular Jaco-type graph was introduced and studied further in [6]. It is well known that for the set $N_0$ of all non-negative integers and $n, k \in N, k \geq 2$ modular arithmetic allows an integer mapping in respect of modulo $k$ as follows:

$$
\begin{aligned}
0 &\mapsto & 0 &= m_0 \\
1 &\mapsto & 1 &= m_1 \\
2 &\mapsto & 2 &= m_2 \\
&\ldots & &\ldots \\
k-1 &\mapsto & k-1 &= m_{k-1} \\
k &\mapsto & k &= m_k \\
k+1 &\mapsto & k+1 &= m_{k+1} \\
&\ldots & &\ldots
\end{aligned}
$$

Note that this new family of Jaco-type graphs, also called the modular Jaco-type graphs, resulting from $\mod k, k \in N$ utilizes a modular non-negative, non-decreasing integer sequence. Let $s_3 = \{a_n\}, a_n \equiv n \pmod{k} = m_n$. Consider the infinite root-graph $J_\infty(s_3)$ and define $d^+(v_i) = m_i$, for $i = 1, 2, 3, \ldots$ Figure 2.3.3 depicts $J_{18}(s_3)$. (Also see [4]).

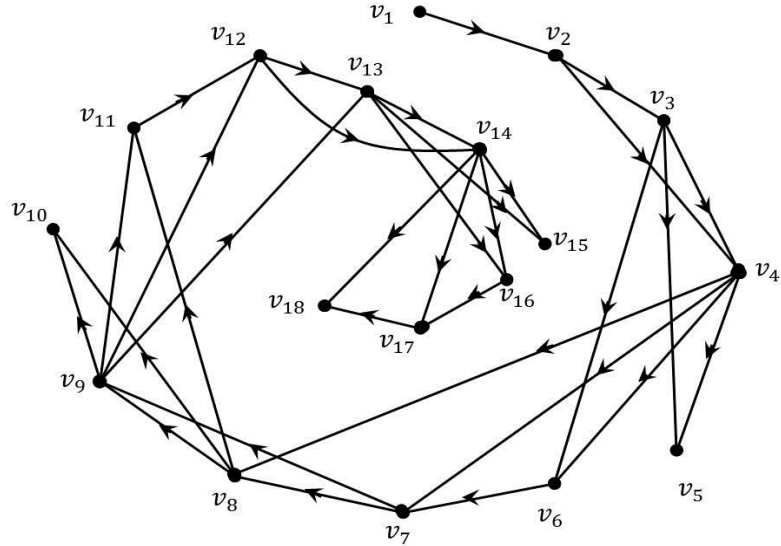

Figure 2.3.3

Table 4. Number of cliques in figure 2.3.3

| $J_n(s_3)$, $n =$ | $\eta^{K_1}(J_n(s_3))$ | $\eta^{K_2}(J_n(s_3))$ | $\eta^{K_3}(J_n(s_3))$ |
|---|---|---|---|
| 1 | 1 | 0 | 0 |
| 2 | 2 | 1 | 0 |
| 3 | 3 | 2 | 0 |
| 4 | 4 | 4 | 1 |
| 5 | 5 | 6 | 2 |
| 6 | 6 | 8 | 3 |
| 7 | 7 | 10 | 4 |
| 8 | 8 | 12 | 5 |
| 9 | 9 | 14 | 6 |
| 10 | 10 | 16 | 7 |
| 11 | 11 | 18 | 8 |
| 12 | 12 | 20 | 9 |
| 13 | 13 | 22 | 10 |
| 14 | 14 | 24 | 11 |
| 15 | 15 | 26 | 12 |
| 16 | 16 | 28 | 13 |
| 17 | 17 | 30 | 14 |
| 18 | 18 | 32 | 15 |

It has been shown in [6] that $d^-(v_n) = 2, n \geq 4$ in $J_n(s_3)$. So on each expansion to the next order $n + 1$, exactly one addition vertex (clique $K_1$) and one addition clique $K_3$ are added. Therefore, $\eta^{K_3}(J_n(s_3)) = 2\eta^{K_3}(J_n(s_3)) + 1, n \geq 3$.

### 2.3.4  Application to the Set Jaco-type Graph: $J_n(s_4)$

We introduce the concept of a set Jaco-type graph. For the set say $A = \{1,2,3,...,n\}$ we have the non-empty subsets by convention in the order $\{1\},\{2\},\{3\},...,\{n\}$, $\{1,2\},\{1,3\},...,\{1,n\}$, $\{2,3\},...,\{2,n\},...$, $\{1,2,3,...,n\}$. Map the vertices, $v_1 \mapsto \{1\}$, $v_2 \mapsto \{2\}$, $v_3 \mapsto \{3\}$, $v_4 \mapsto \{1,2\}$, $v_5 \mapsto \{1,3\},...$, $v_{2^n-1} \mapsto \{1,2,3,...,v_n\}$. Define the out-degree of vertex $v_i$ to be sum of elements of subset $i$. Let the vertex degree of vertex $v_j, j > 2^n - 1$ mapped onto corresponding vertex degree of vertex subscript $i = 1 + (j-1) \bmod (2^n - 1)$.

This graph is called the *set Jaco-type graph.* Figure 2.3.4 depicts the set Jaco-type graph $J_{13}(s_4)$ with the terms of the sequence $\{s_4\}$ corresponding to $a_1 \mapsto 1$, $a_2 \mapsto 2$, $a_3 \mapsto 3$, $a_4 \mapsto 4$, $a_5 \mapsto 4$, $a_6 \mapsto 5$, $a_7 \mapsto 6$, $a_8 \mapsto 1$, $a_9 \mapsto 2$, $\cdots$, $a_i \mapsto 1 + (j-1) \bmod (2^n - 1)$.

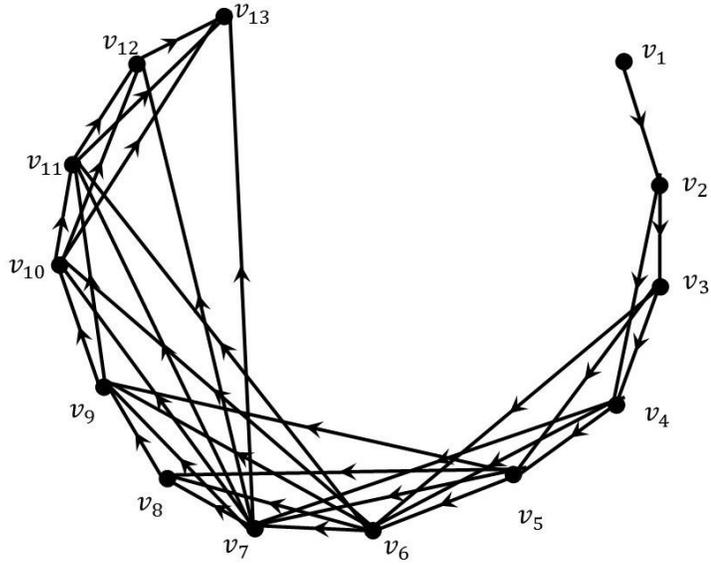

Figure 2.3.4

Table 5. Number of cliques found in figure 2.3.4

| $K_n$, $n=$ | $\eta^{K_1}(J_n(s_4))$ | $\eta^{K_2}(J_n(s_4))$ | $\eta^{K_3}(J_n(s_4))$ | $\eta^{K_4}(J_n(s_4))$ | $\eta^{K_5}(J_n(s_4))$ |
|---|---|---|---|---|---|
| 1 | 1 | 0 | 0 | 0 | 0 |
| 2 | 2 | 1 | 0 | 0 | 0 |
| 3 | 3 | 2 | 0 | 0 | 0 |
| 4 | 4 | 3 | 1 | 0 | 0 |
| 5 | 5 | 5 | 2 | 0 | 0 |
| 6 | 6 | 8 | 5 | 1 | 0 |
| 7 | 7 | 11 | 8 | 2 | 0 |
| 8 | 8 | 14 | 11 | 3 | 0 |
| 9 | 9 | 18 | 17 | 7 | 1 |
| 10 | 10 | 21 | 20 | 8 | 1 |
| 11 | 11 | 25 | 26 | 12 | 2 |
| 12 | 12 | 28 | 29 | 13 | 2 |
| 13 | 13 | 32 | 36 | 17 | 3 |

Note that because the vertex degrees are specifically defined as the sum of the elements of proper subsets which is a function of the cardinality of a given set, it is not possible to determine a closed formula for $\eta^{K_i}(J_n(s_4))$ in general. The modular Jaco-type graph has a similar limitation on generality.

## 3. CONCLUSION

This paper discusses the introduction to the concept of Jaco-type graphs with clique parameters. It includes theorems to find *l*-cliques in certain Jaco-type graphs. Applications to certain Jaco-type graphs in respect of the concept of vertex clique degrees were also presented.